\documentclass[11pt, letterpaper, reqno]{article}

\usepackage[utf8]{inputenc}
\usepackage[T1]{fontenc}
\usepackage{lmodern}
\usepackage[final]{microtype}

\usepackage{amsmath,amsthm,amssymb,amsfonts}
\usepackage{mathtools}
\usepackage{bm}

\usepackage{graphicx}
\usepackage{float}
\usepackage{booktabs}
\usepackage{caption}
\usepackage{subcaption}

\usepackage[margin=1.1in]{geometry}
\usepackage[shortlabels]{enumitem}
\usepackage[dvipsnames]{xcolor}

\usepackage{hyperref}
\usepackage[nameinlink, capitalize]{cleveref}

\theoremstyle{plain}
\newtheorem{theorem}{Theorem}[section]

\newtheorem{conjecture}[theorem]{Conjecture}

\theoremstyle{definition}
\newtheorem{definition}[theorem]{Definition}

\theoremstyle{remark}

\DeclareMathOperator{\GP}{GP}
\DeclareMathOperator{\Ind}{I}
\DeclareMathOperator{\Tr}{Tr}

\hypersetup{
    colorlinks=true,
    linkcolor=NavyBlue,
    citecolor=ForestGreen,
    urlcolor=NavyBlue,
    pdftitle={Parity-Dependent Real-Rootedness in Independence Polynomials},
    pdfauthor={Rohan Pandey}
}

\title{\textbf{Parity-Dependent Real-Rootedness in Independence Polynomials of Generalized Petersen Graphs}}

\author{
    Rohan Pandey \\
    \small University of Washington \\
    \small \texttt{rpande@uw.edu}
}
\date{\today}

\begin{document}

\maketitle

\begin{abstract}
We investigate the distribution of zeros of the independence polynomial $\Ind(G, x)$ for the family of Generalized Petersen graphs $\GP(n, k)$ in the complex plane. While the independence numbers and coefficients of these graphs have been studied, the global behavior of their roots remains largely unexplored. Using an exact transfer matrix algorithm parameterized by $k$, we compute $\Ind(\GP(n,k), x)$ for $n$ up to $30$ and $k \in \{1, 2, 3, 4\}$. Our numerical analysis reveals a striking parity-based dichotomy: for odd $k$, the roots exhibit complex conjugate structures accumulating on closed curves, whereas for even $k$, the roots appear to be strictly real and negative. Motivated by this evidence, we conjecture that $\Ind(\GP(n,k), x)$ is real-rooted, and hence log-concave, if and only if $k$ is even. This phenomenon connects algebraic properties of $\GP(n,k)$ to questions about zero-free regions and limiting behavior in the hard-core lattice gas model.
  

\end{abstract}

\section{Introduction}

The independence polynomial of a graph encodes fundamental combinatorial information about its independent sets. Given a graph $G$, the \emph{independence polynomial} is defined as
\begin{equation}
    \Ind(G, x) = \sum_{k=0}^{\alpha(G)} i_k(G) \, x^k,
\end{equation}
where $i_k(G)$ denotes the number of independent sets of cardinality $k$ in $G$, and $\alpha(G)$ is the independence number~\cite{gutman1983independence}.

The study of these polynomials is motivated by applications in statistical mechanics, where $\Ind(G, x)$ coincides with the partition function of the hard-core lattice gas model with fugacity $x$~\cite{scott2005repulsive}. A central problem is determining the location of the polynomial’s zeros in the complex plane. In physics, these zeros (Lee–Yang zeros) govern phase transitions: the absence of roots in a neighborhood of the positive real axis implies analyticity of the free energy~\cite{heilmann1972theory}. In combinatorics, if all roots of $\Ind(G,x)$ are real, Newton’s inequalities imply that the coefficient sequence $\{i_k(G)\}$ is log-concave and unimodal~\cite{stanley1989log}.

Real-rootedness has been established for matching polynomials of line graphs~\cite{heilmann1972theory} and for independence polynomials of claw-free graphs~\cite{chudnovsky2007roots}. In contrast, much less is known for cubic graphs that fall outside these classes. The Generalized Petersen graphs $\GP(n,k)$ form a classical family of cubic graphs. While their independence numbers have been determined~\cite{fox2012independence,besharati2010independence}, the global structure of their independence polynomial roots has not previously been studied.

In this work, we conduct a computational investigation of the roots of $\Ind(\GP(n,k),x)$ using an exact transfer matrix framework. Our results reveal a robust parity phenomenon: the polynomials appear to be real-rooted precisely when the step size $k$ is even. Since our conclusions are based on finite but exact computations and numerical root-finding, they motivate a conjecture rather than a proof. The main contribution is the identification of a previously unreported parity dichotomy in the root geometry of $\Ind(\GP(n,k),x)$, together with a transfer-matrix framework that enables exact computation in this family.

\paragraph{Scope and nature of results.}
This work is computational in nature, based on exact transfer-matrix constructions and numerical root-finding for finite but systematically increasing values of $n$. While no proofs are given, the consistency of the observed behavior across all tested instances motivates a precise conjecture. The goal of this paper is to identify and clearly formulate this phenomenon as a potential target for future theoretical investigation.

\section{Preliminaries and Methodology}

\subsection{Generalized Petersen Graphs}

\begin{definition}[\cite{watkins1969theorem}]
For integers $n \ge 3$ and $1 \le k < n/2$, the graph $\GP(n,k)$ has vertex set $U \cup V$, where
\[
U = \{u_0,\dots,u_{n-1}\}, \quad V = \{v_0,\dots,v_{n-1}\}.
\]
The edge set consists of
\begin{enumerate}[(i)]
    \item $u_i u_{i+1}$ (outer cycle),
    \item $v_i v_{i+k}$ (inner chords),
    \item $u_i v_i$ (spokes),
\end{enumerate}
with indices taken modulo $n$.
\end{definition}

\subsection{Parameterized Transfer Matrix Method}

To compute $\Ind(\GP(n,k),x)$ for large $n$, we generalize the transfer matrix method originally developed for ladder and grid graphs~\cite{calkin2003counting}. Unlike the prism graph ($k=1$), the presence of long-range inner edges in $\GP(n,k)$ requires tracking the occupancy history of multiple vertices.

At step $j$, we encode the state as a binary vector
\[
S_j = (u_j, v_j, v_{j-1}, \dots, v_{j-k+1}) \in \{0,1\}^{k+1},
\]
where each entry indicates whether the corresponding vertex is included in the independent set. The transfer matrix $T_k$ has dimension $2^{k+1} \times 2^{k+1}$. An entry $T_k(A,B)$ is nonzero, with weight $x^{|B|}$, if and only if state $B$ can validly follow state $A$, subject to:
\begin{itemize}
    \item \textbf{Shift consistency:} The history components of $B$ agree with the trailing components of $A$.
    \item \textbf{Independence constraints:} No edge $(u_j,u_{j-1})$, $(u_j,v_j)$, or $(v_j,v_{j-k})$ has both endpoints occupied.
\end{itemize}

Cyclic closure is enforced by taking the trace of $T_k^n$, yielding
\begin{equation}
    \Ind(\GP(n,k),x) = \Tr(T_k^n).
\end{equation}
All coefficients are computed exactly as integers. Roots are then obtained via standard polynomial root-finding, with numerical precision verified to $10^{-10}$.

\paragraph{Computational regime and verification.}
All transfer matrices were constructed exactly, and polynomial coefficients were computed symbolically as integers. Roots were obtained using standard polynomial root-finding routines and verified to numerical precision $10^{-10}$. The computations were carried out for all $k \in \{1,2,3,4\}$ and for $20 \le n \le 30$, with consistent qualitative behavior observed throughout this range.

\section{Root Distribution Analysis}

\begin{figure}[htbp]
    \centering
    \begin{subfigure}[b]{0.48\textwidth}
        \centering
        \includegraphics[width=\linewidth]{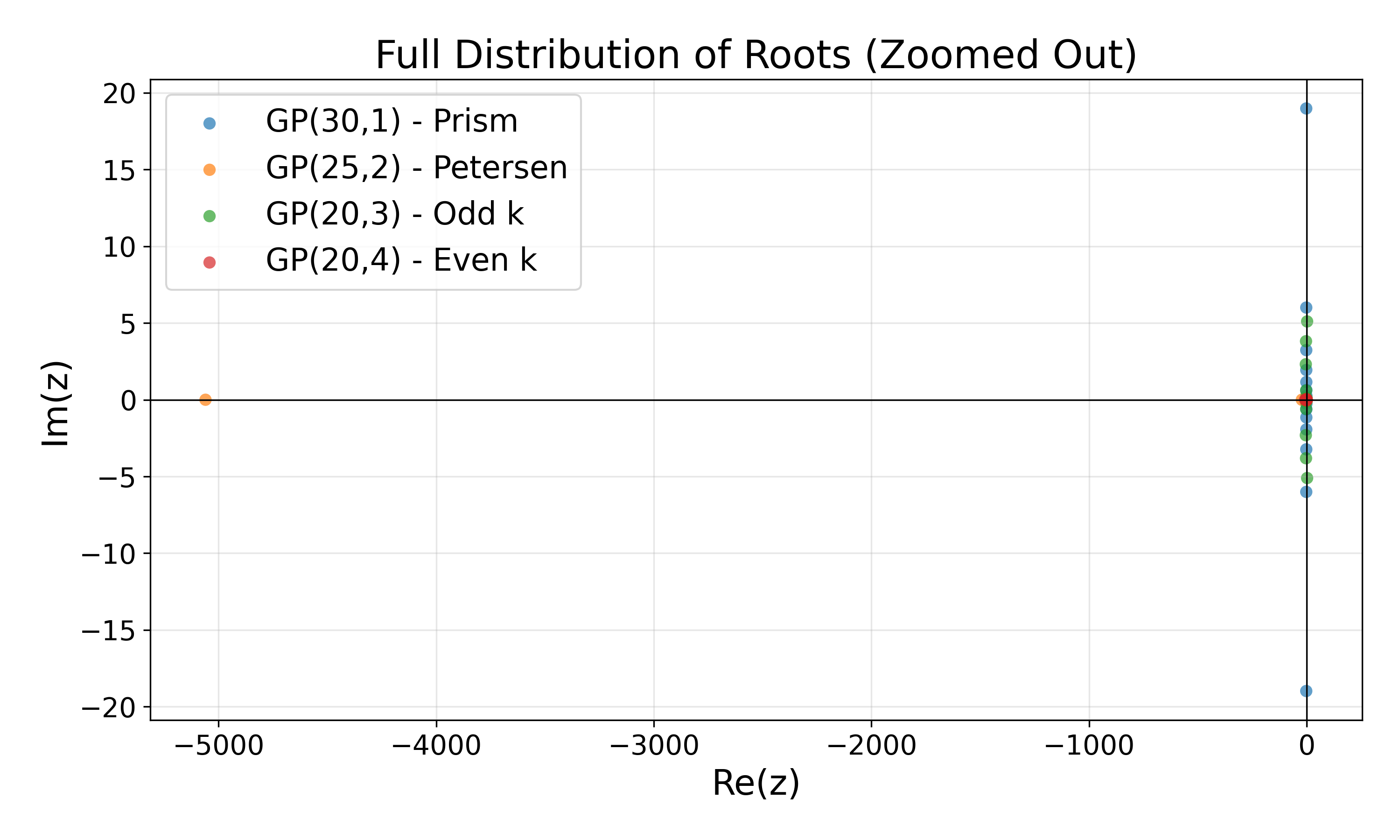}
        \caption{Global view showing real roots extending along the negative real axis.}
    \end{subfigure}
    \hfill
    \begin{subfigure}[b]{0.48\textwidth}
        \centering
        \includegraphics[width=\linewidth]{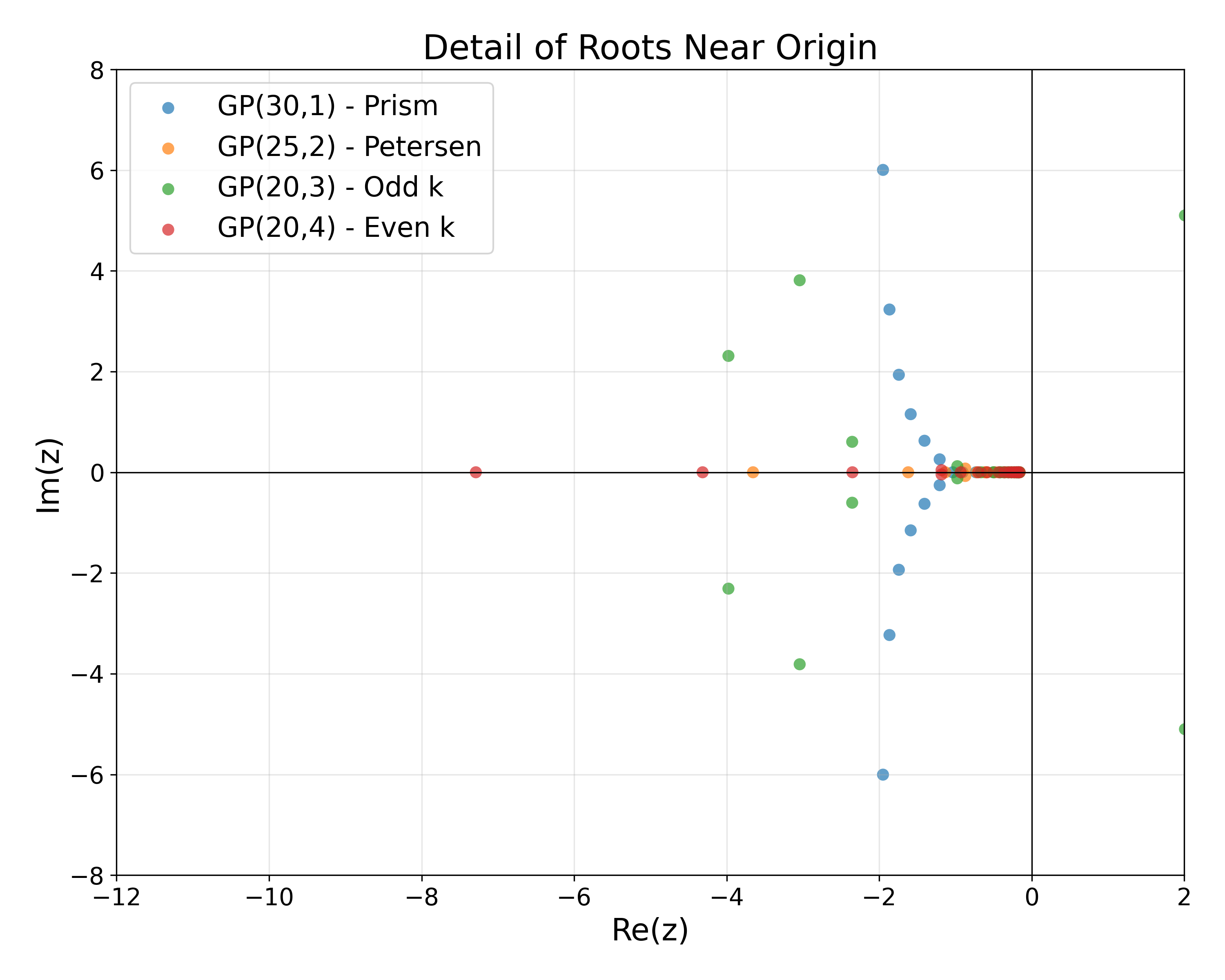}
        \caption{Zoom near the origin highlighting complex conjugate structures for odd $k$.}
    \end{subfigure}
    \caption{\textbf{Root distributions for $\GP(n,k)$.} Even values of $k$ produce exclusively real negative roots, while odd values of $k$ generate complex conjugate curves.}
    \label{fig:roots}
\end{figure}

We computed the roots of $\Ind(\GP(n,k),x)$ across the computational regime described in Section~2. Representative results are shown in \Cref{fig:roots}.

\subsection{Odd $k$: Complex Roots and Limit Shapes}

For $k=1$ and $k=3$, the roots consistently appear off the real axis.
\begin{itemize}
    \item \textbf{$k=1$:} Roots form a closed vertical oval centered near $\Re(z) \approx -0.5$, confirming that $\GP(n,1)$ is not real-rooted.
    \item \textbf{$k=3$:} Roots form more intricate conjugate curves with a distinct lobe structure.
\end{itemize}
As $n$ increases, these roots appear to converge toward continuous limiting curves, consistent with phenomena observed for other graph families~\cite{brown2007roots}.

\subsection{Even $k$: Apparent Real-Rootedness}

In contrast, for $k=2$ and $k=4$, all computed roots lie on the negative real axis.
\begin{itemize}
    \item \textbf{$k=2$:} For $n=25$, all $2n$ roots were real within numerical tolerance.
    \item \textbf{$k=4$:} For $n=20$, the same behavior was observed.
\end{itemize}
This persistence strongly suggests genuine real-rootedness rather than a finite-size artifact.

\section{Discussion and Conjecture}

The parity-dependent behavior likely arises from structural differences in the inner edge configuration. When $k$ is even, the inner subgraph decomposes into bipartite components, while odd $k$ introduces odd cycles that may obstruct real-rootedness. Although existing techniques such as interlacing or stability arguments do not directly apply, the sharpness of the observed dichotomy suggests that a structural explanation may be accessible.

\begin{conjecture}[Parity Conjecture]
For all integers $n \ge 2k+1$, the independence polynomial $\Ind(\GP(n,k),x)$ has only real roots if and only if $k$ is even.
\end{conjecture}

\paragraph{Significance.}
If proven, this conjecture would identify $\GP(n,2k)$ as a new infinite family of graphs with real-rooted independence polynomials outside the claw-free class, helping to clarify the boundary between real-rooted and non-real-rooted behavior for cubic graphs.

\section{Conclusion}

We have identified a parity-driven dichotomy in the root distributions of independence polynomials for Generalized Petersen graphs. Odd step sizes yield complex root curves, while even step sizes appear to enforce real-rootedness. These findings motivate further theoretical investigation into the algebraic mechanisms underlying this phenomenon.

\section*{Data Availability}
The source code used to generate all transfer matrices and root plots is publicly available at  
\url{https://github.com/Rohan-Pandey1729/polynomial-independence}.

\section*{Acknowledgements}
The author thanks Professor Konstantinos Mamis for his guidance and Ray Chen for helpful discussions.

\bibliographystyle{plain}
\bibliography{references}

\end{document}